\documentclass[12pt]{article}
\usepackage{amsmath}
\usepackage{amsfonts}
\usepackage{mitpress}
\usepackage{graphicx}

\begin{document}

\title{Galerkin Method for the numerical solution of the advection-diffusion
equation by using exponential B-splines}
\author{Melis Zorsahin Gorgulu and Idris Dag \\
{\small Department of Mathematics-Computer, Eskisehir Osmangazi University,
26480, Eskisehir, Turkey.}}
\maketitle

\begin{abstract}
In this paper, the exponential B-spline functions are used for the numerical
solution of the advection-diffusion equation. Two numerical examples\
related to pure advection in a finitely long channel and the distribution of
an initial Gaussian pulse are employed to illustrate the accuracy and the
efficiency of the method. Obtained results are compared with some early
studies.

\end{abstract}

\section{Introduction}

One known problem of our age is the environmental pollution. This problem is
increasingly reduce the quality of our water. Scientists are benefiting from
the solution of the advection-diffusion equation (ADE) in modelling this
problem. Since both advection and diffusion terms exist in the ADE, it also
arises very frequently in transferring mass, heat, energy, velocity and
vorticity in engineering and chemistry. Thus, the heat transfer in a
draining film, dispersion of tracers in porous media, the intrusion of salt
water into fresh water aquifers, the spread of pollutants in rivers and
streams, the dispersion of dissolved material in estuaries and coastal sea,
contaminant dispersion in shallow lakes, the absorption of chemicals into
beds, the spread of solute in a liquid flowing through a tube, long-range
transport of pollutants in the atmosphere, forced cooling by fluids of solid
material such as windings in turbo generators, thermal pollution in river
systems and flow in porous media, etc. are modeled by ADE \cite{korkmaz}.

It is well known that the solution of the advection-diffusion boundary value
problem displays sharp boundary layers. To cope with the sharp solutions,
some of the spline based methods for the numerical solution of ADE are
suggested such as the quasi-Lagrangian cubic spline method \cite%
{pepper,okamoto}, the characteristic methods integrated with splines \cite%
{szym,tsai}, the cubic B-spline Galerkin method \cite{gardner}, the
quadratic B-spline subdomain collocation method \cite{gardner3}, the spline
approximation with the help of upwind collocation nodes \cite{funaro}, the
exponential spline interpolation in characteristic based scheme \cite{zoppou}%
, the cubic spline interpolation for the advection component and the
Crank-Nicolson scheme for the diffusion component \cite{ahmad,koth}, the
meshless method based on thin-plate spline radial basis functions \cite%
{boztosun}, the least-square B-spline finite element method \cite{dag}, the
standard finite difference method \cite{thong,thong2}, the cubic B-spline
collocation method \cite{goh,goh2,dirk}, the quadratic/cubic\ B-spline
Taylor-Galerkin methods \cite{dag2}, the cubic/quadratic B-spline
least-squares finite element techniques \cite{kapoor,dhawan,dhawan2}, the
cubic B-spline differential quadrature method \cite{korkmaz}, the quadratic
Galerkin finite elements method \cite{bulut}.

The exponential B-spline basis functions are used to establish the numerical
methods. Thus the exponential B-spline based collocation method are
constructed to solve the differential equations. Numerical solution of the
singular perturbation problem is solved with a variant of exponential
B-spline collocation method in the work \cite{ms1}, the cardinal exponential
B-splines is used for solving the singularly perturbed problems \cite{ca},
the exponential B-spline collocation method is built up for finding the
numerical solutions of the self-adjoint singularly perturbed boundary value
problems in the work \cite{rao}, the numerical solutions of the
convection-diffusion equation is obtained by using the exponential B-spline
collocation method \cite{mohammadi}. As far as we search, no study exists
solving the advection-diffusion problems using the exponential B-spline
Galerkin method. Thus advection-diffusion equation is fully integrated with
combination of the exponential B-spline Galerkin method (EBSGM) for space
discretization and Crank-Nicolson method for time discretization.

The study is organized as follows. In section 2, exponential B-splines are
introduced and their some basic relations are given. In section 3, the
application of the numerical method to the ADE is given. The efficiency and
the accuracy of the present method are investigated by using two numerical
experiments related to pure advection in an infinitely long channel and the
distribution of an initial Gaussian pulse.

\section{Exponential B-splines and Finite Element Solution}

The mathematical model describing the transport and diffusion processes is
the one dimensional ADE%
\begin{equation}
\frac{\partial u}{\partial t}+\xi \frac{\partial u}{\partial x}-\lambda
\frac{\partial ^{2}u}{\partial x^{2}}=0,  \label{1}
\end{equation}%
where the function $u(x,t)$ represents the concentration at position $x$ and
time $t$ with uniform flow velocity $\xi $ and constant diffusion
coefficient $\lambda .$ The initial condition of Eq. (\ref{1}) is%
\begin{equation}
u(x,0)=u_{0}(x),\text{ \ }0\leq x\leq L  \label{2a}
\end{equation}%
and the boundary conditions are%
\begin{equation}
u(0,t)=f_{0}(t),\text{ \  \ }u(L,t)=f_{L}(t)\text{ \  \ or }-\lambda \dfrac{%
\partial u}{\partial x}|_{L}=\phi _{L}(t)  \label{2b}
\end{equation}%
where $L$ is the length of the channel, $\phi _{L}$ is the flux at the
boundary $x=L$ and $u_{0},$ $f_{0},$ $f_{L}$ are imposed functions.

Let us consider a uniform mesh $\Gamma $ with the knots $x_{i}$ on $\left[
a,b\right] $ such that%
\begin{equation*}
\Gamma :a=x_{0}<x_{1}<x_{2}<\cdots <x_{N-1}<x_{N}=b
\end{equation*}%
\noindent where $h=\dfrac{b-a}{N}$ and $x_{i}=x_{0}+ih.$

Let $\phi _{i}\left( x\right) $ be the exponential B-splines at the points
of $\Gamma $ together with knots $x_{i}$, $i=-3,-2,-1,N+1,N+2,N+3$ outside
the interval $[a,b]$ and having a finite support on the four consecutive
intervals $\left[ x_{i}+kh,x_{i}+\left( k+1\right) h\right] _{k=-3}^{0},$ $%
i=0,...,N+2.$ According to McCartin \cite{mccartin}, the $\phi _{i}\left(
x\right) $ can be defined as

\begin{equation}
\phi _{i}\left( x\right) =\left \{
\begin{array}{lll}
b_{2}\left[ \left( x_{i-2}-x\right) -\dfrac{1}{p}\left( \sinh \left( p\left(
x_{i-2}-x\right) \right) \right) \right] & \text{ \ } & \text{if }x\in \left[
x_{i-2},x_{i-1}\right] ; \\
a_{1}+b_{1}\left( x_{i}-x\right) +c_{1}e^{p\left( x_{i}-x\right)
}+d_{1}e^{-p\left( x_{i}-x\right) } & \text{ } & \text{if }x\in \left[
x_{i-1},x_{i}\right] ; \\
a_{1}+b_{1}\left( x-x_{i}\right) +c_{1}e^{p\left( x-x_{i}\right)
}+d_{1}e^{-p\left( x-x_{i}\right) } & \text{ } & \text{if }x\in \left[
x_{i},x_{i+1}\right] ; \\
b_{2}\left[ \left( x-x_{i+2}\right) -\dfrac{1}{p}\left( \sinh \left( p\left(
x-x_{i+2}\right) \right) \right) \right] & \text{ } & \text{if }x\in \left[
x_{i+1},x_{i+2}\right] ; \\
0 & \text{ } & \text{otherwise.}%
\end{array}%
\right.  \label{3}
\end{equation}%
where

\begin{equation*}
\begin{array}{l}
p=\underset{0\leq i\leq N}{\max }p_{i},\text{ }s=\sinh \left( ph\right) ,%
\text{ }c=\cosh \left( ph\right) \\
b_{2}=\dfrac{p}{2\left( phc-s\right) },\text{ }a_{1}=\dfrac{phc}{phc-s},%
\text{ }b_{1}=\dfrac{p}{2}\left[ \dfrac{c\left( c-1\right) +s^{2}}{\left(
phc-s\right) \left( 1-c\right) }\right] , \\
c_{1}=\dfrac{1}{4}\left[ \dfrac{e^{-ph}\left( 1-c\right) +s\left(
e^{-ph}-1\right) }{\left( phc-s\right) \left( 1-c\right) }\right] ,\text{ }%
d_{1}=\dfrac{1}{4}\left[ \dfrac{e^{ph}\left( c-1\right) +s\left(
e^{ph}-1\right) }{\left( phc-s\right) \left( 1-c\right) }\right] .%
\end{array}%
\end{equation*}

Each basis function $\phi _{i}\left( x\right) $ is twice continuously
differentiable. The values of $\phi _{i}\left( x\right) ,$ $\phi
_{i}^{\prime }\left( x\right) $ and $\phi _{i}^{\prime \prime }\left(
x\right) $ at the knots $x_{i}$'s are given in Table 1.%
\begin{equation*}
\begin{tabular}{c|ccccc}
\multicolumn{6}{l}{Table 1: Exponential B-spline values} \\ \hline \hline
& $x_{i-2}$ & $x_{i-1}$ & $x_{i}$ & $x_{i+1}$ & $x_{i+2}$ \\ \hline
$\phi _{i}\left( x\right) $ & $0$ & $\frac{s-ph}{2\left( phc-s\right) }$ & $%
1 $ & $\frac{s-ph}{2\left( phc-s\right) }$ & $0$ \\
$\phi _{i}^{\prime }\left( x\right) $ & $0$ & $\frac{p\left( c-1\right) }{%
2\left( phc-s\right) }$ & $0$ & $\frac{p\left( 1-c\right) }{2\left(
phc-s\right) }$ & $0$ \\
$\phi _{i}^{\prime \prime }\left( x\right) $ & $0$ & $\frac{p^{2}s}{2\left(
phc-s\right) }$ & $\frac{-p^{2}s}{phc-s}$ & $\frac{p^{2}s}{2\left(
phc-s\right) }$ & $0$ \\ \hline \hline
\end{tabular}%
\end{equation*}

The $\phi _{i}\left( x\right) ,$ $i=-1,\ldots ,N+1$ form a basis for
functions defined on the interval $[a,b]$. We seek an approximation $U(x,t)$
to the analytical solution $u(x,t)$ in terms for the exponential B-splines%
\begin{equation}
u\left( x,t\right) \approx U\left( x,t\right) =\overset{N+1}{\underset{i=-1}{%
\sum }}\phi _{i}\left( x\right) \delta _{i}\left( t\right)  \label{4}
\end{equation}%
where $\delta _{i}\left( t\right) $ are time dependent unknown to be
determined from the boundary conditions and Galerkin approach to the
equation (\ref{1}). The approximate solution and their derivatives at the
knots can be found from the Eq. (\ref{3}-\ref{4}) as
\begin{equation}
\begin{tabular}{l}
$U_{i}=U(x_{i},t)=\alpha _{1}\delta _{i-1}+\delta _{i}+\alpha _{1}\delta
_{i+1},$ \\
$U_{i}^{\prime }=U^{\prime }(x_{i},t)=\alpha _{2}\delta _{i-1}-\alpha
_{2}\delta _{i+1},$ \\
$U_{i}^{\prime \prime }=U^{\prime \prime }(x_{i},t)=\alpha _{3}\delta
_{i-1}-2\alpha _{3}\delta _{i}+\alpha _{3}\delta _{i+1}$%
\end{tabular}
\label{5}
\end{equation}%
where $\alpha _{1}=\dfrac{s-ph}{2(phc-s)},\alpha _{2}=\dfrac{p(1-c)}{2(phc-s)%
},\alpha _{3}=\dfrac{p^{2}s}{2(phc-s)}.$

Applying the Galerkin method to the ADE with the exponential B-splines as
weight function over the interval $\left[ a,b\right] $ gives
\begin{equation}
\underset{a}{\overset{b}{\int }}\phi _{i}\left( x\right) \left( u_{t}+\xi
u_{x}-\lambda u_{xx}\right) dx=0.  \label{6}
\end{equation}

The approximate solution $U$ over the element $[x_{m},x_{m+1}]$ can be
written as

\begin{equation}
U^{e}=\phi _{m-1}\left( x\right) \delta _{m-1}\left( t\right) +\phi
_{m}\left( x\right) \delta _{m}\left( t\right) +\phi _{m+1}\left( x\right)
\delta _{m+1}\left( t\right) +\phi _{m+2}\left( x\right) \delta _{m+2}\left(
t\right)  \label{7}
\end{equation}%
where quantities $\delta _{j}\left( t\right) ,$ $j=m-1,...,m+2$ are element
parameters and $\phi _{j}\left( x\right) ,$ $j=m-1,...,m+2$ are known as the
element shape functions.

The contribution of the integral equation (\ref{6}) over the sample interval
$[x_{m},x_{m+1}]$ is given by
\begin{equation}
\underset{x_{m}}{\overset{x_{m+1}}{\int }}\phi _{j}\left( x\right) \left(
u_{t}+\xi u_{x}-\lambda u_{xx}\right) dx.  \label{8}
\end{equation}%
Applying the Galerkin discretization scheme by replacing $U_{t},$ $U_{x},$ $%
U_{xx},$ which are derivatives of the approximate solution $U^{e}$ in Eq. (%
\ref{7}), into $u_{t},$ $u_{x},$ $u_{xx},$ which are derivatives of the
exact solution $u$, respectively, we obtain a system of equations in the
unknown parameters $\delta _{j}$
\begin{equation}
\overset{m+2}{\underset{i=m-1}{\sum }}\left \{ \left( \underset{x_{m}}{%
\overset{x_{m+1}}{\int }}\phi _{j}\phi _{i}dx\right) \overset{\mathbf{%
\bullet }}{\delta }_{i}+\xi \left( \underset{x_{m}}{\overset{x_{m+1}}{\int }%
}\phi _{j}\phi _{i}^{\prime }dx\right) \delta _{i}-\lambda \left( \underset{%
x_{m}}{\overset{x_{m+1}}{\int }}\phi _{j}\phi _{i}^{\prime \prime
}dx\right) \delta _{i}\right \}  \label{9}
\end{equation}%
where $i$ and $j$ take only the values $m-1,$ $m,$ $m+1,$ $m+2$ for $%
m=0,1,\ldots ,N-1$ and $\overset{\mathbf{\bullet }}{}$ denotes time
derivative.

In the above system of differential equations, when $A_{ji}^{e},$ $%
B_{ji}^{e} $ and $C_{ji}^{e}$ are denoted by%
\begin{equation}
\begin{tabular}{lll}
$A_{ji}^{e}=\underset{x_{m}}{\overset{x_{m+1}}{\int }}\phi _{j}\phi _{i}dx,$
& $B_{ji}^{e}=\underset{x_{m}}{\overset{x_{m+1}}{\int }}\phi _{j}\phi
_{i}^{\prime }dx,$ & $C_{ji}^{e}=\underset{x_{m}}{\overset{x_{m+1}}{\int }}%
\phi _{j}\phi _{i}^{\prime \prime }dx$%
\end{tabular}
\label{10}
\end{equation}%
where $\mathbf{A}^{e},$ $\mathbf{B}^{e}$ and $\mathbf{C}^{e}$ are the
element matrices of which dimensions are $4\times 4$, the matrix form of the
Eq.(\ref{9}) can be written as%
\begin{equation}
\mathbf{A}^{e}\overset{\mathbf{\bullet }}{\mathbf{\delta }^{e}}+\left( \xi
\mathbf{B}^{e}-\lambda \mathbf{C}^{e}\right) \mathbf{\delta }^{e}  \label{11}
\end{equation}%
where $\mathbf{\delta }^{e}\mathbf{=}\left( \delta _{m-1},...,\delta
_{m+2}\right) ^{T}.$

Gathering the systems (\ref{11}) over all elements, we obtain global system%
\begin{equation}
\mathbf{A}\overset{\mathbf{\bullet }}{\mathbf{\delta }}+\left( \xi \mathbf{B}%
-\lambda \mathbf{C}\right) \mathbf{\delta }=0  \label{12}
\end{equation}%
where $\mathbf{A},$ $\mathbf{B},$ $\mathbf{C}$ are derived from the
corresponding element matrices $\mathbf{A}^{e},$ $\mathbf{B}^{e},$ $\mathbf{C%
}^{e}$ and $\mathbf{\delta =}\left( \delta _{-1},...,\delta _{N+1}\right)
^{T}$ contain all elements parameters.

The unknown parameters $\mathbf{\delta }$ are interpolated between two time
levels $n$ and $n+1$ with the Crank-Nicolson method%
\begin{equation*}
\begin{array}{cc}
\mathbf{\delta }=\dfrac{\delta ^{n+1}+\delta ^{n}}{2}, & \overset{\mathbf{%
\bullet }}{\mathbf{\delta }}=\dfrac{\delta ^{n+1}-\delta ^{n}}{\Delta t},%
\end{array}%
\end{equation*}%
we obtain iterative formula for the time parameters $\mathbf{\delta }^{n}$:%
\begin{equation}
\left[ \mathbf{A+}\frac{\Delta t}{2}\left( \xi \mathbf{B}-\lambda \mathbf{C}%
\right) \right] \mathbf{\delta }^{n+1}=\left[ \mathbf{A-}\frac{\Delta t}{2}%
\left( \xi \mathbf{B}-\lambda \mathbf{C}\right) \right] \mathbf{\delta }^{n}.
\label{13}
\end{equation}%
The set of equations consist of $\left( N+3\right) $ equations with $\left(
N+3\right) $ unknown parameters. Before starting the iteration procedure,
boundary conditions must be adapted into the system. For this purpose, we
delete first and last equations from the system (\ref{13}) and eliminate the
terms $\delta _{-1}^{n+1}$ and $\delta _{N+1}^{n+1}$ from the remaining
system (\ref{13}) by using boundary conditions in (\ref{2b}), which give the
following equations:%
\begin{equation*}
\begin{array}{l}
u\left( a,t\right) =\alpha _{1}\delta _{-1}^{n}+\delta _{0}^{n}+\alpha
_{1}\delta _{1}^{n}=\beta _{1}, \\
u\left( b,t\right) =\alpha _{1}\delta _{N-1}^{n}+\delta _{N}^{n}+\alpha
_{1}\delta _{N+1}^{n}=\beta _{2},%
\end{array}%
\end{equation*}%
we obtain a septa-diagonal matrix with the dimension $\left( N+1\right)
\times \left( N+1\right) $.

To start evolution of the vector of initial parameters $\mathbf{\delta }^{0}$%
, it must be determined by using the initial condition (\ref{2a}) and
boundary conditions (\ref{2b}):%
\begin{equation}
\begin{tabular}{l}
$u_{0}^{\prime }(x_{0},0)=\dfrac{p\left( 1-c\right) }{2\left( phc-s\right) }%
\delta _{-1}+\dfrac{p\left( c-1\right) }{2\left( phc-s\right) }\delta _{1}$
\\
$u\left( x_{m},0\right) =\dfrac{s-ph}{2\left( phc-s\right) }\delta
_{m-1}+\delta _{m}+\dfrac{s-ph}{2\left( phc-s\right) }\delta _{m+1},$ $%
m=0,...,N$ \\
$u^{\prime }\left( x_{N},0\right) =\dfrac{p\left( 1-c\right) }{2\left(
phc-s\right) }\delta _{N-1}+\dfrac{p\left( c-1\right) }{2\left( phc-s\right)
}\delta _{N+1}$%
\end{tabular}
\label{15}
\end{equation}%
The solution of matrix equation (\ref{15}) with the dimensions $\left(
N+1\right) \times \left( N+1\right) $ is obtained by the way of Thomas
algorithm. Once $\mathbf{\delta }^{0}$ is determined, we can start the
iteration of the system to find the parameters $\mathbf{\delta }^{n}$ at
time $t^{n}=n\Delta t.$ Thus the approximate solution $U$ (\ref{4}) can be
determined by using these $\mathbf{\delta }$ values.

\section{Test Problems}

We have carried out two test problems to demonstrate the performance of the
given algorithm. Accuracy of the method is measured by the error norm%
\begin{equation}
L_{\infty }=\left \Vert u^{\text{exact}}-u^{\text{numeric}}\right \Vert
_{\infty }=\max_{0\leq j\leq N}\left \vert u_{j}^{\text{exact}}-u_{j}^{\text{%
numeric}}\right \vert .  \label{15a}
\end{equation}%
In numerical calculations, the determination of $p$ in the exponential
B-spline is made by experimentally. The Courant number is defined by the
ratio of the flow velocity $\xi $ to the mesh velocity $\dfrac{h}{\Delta t},$
i.e.,
\begin{equation*}
C_{r}=\xi \dfrac{\Delta t}{h}.
\end{equation*}

\subsection{Pure advection in an infinitely long channel}

\begin{figure}[ht]
\centering\includegraphics[scale=0.4]{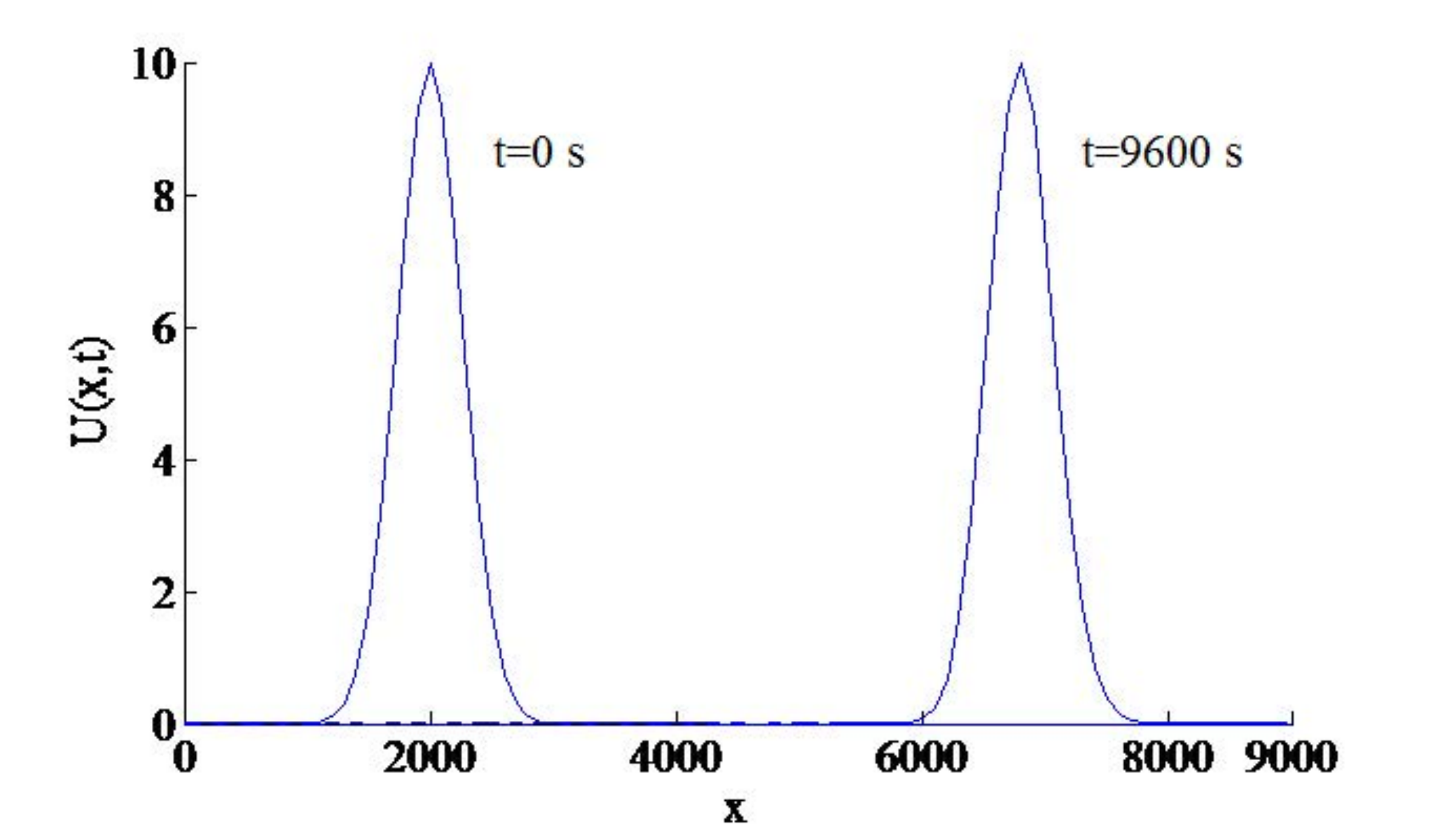}
\caption{{}{\protect \ Transportation of the initial distribution with $C_{r}=0.25$ and $\Delta t=50$}}
\label{fig1}
\end{figure}

In the first example, we consider the pure advection that is $\lambda =0,$
in an infinitely long channel is of long constant cross-section, bottom
slope and in which constant velocity is $\xi =0.5$ $m/s$. The analytical
solution is%
\begin{equation}
u(x,t)=10\exp \left( -\frac{1}{2\rho ^{2}}\left( x-x_{0}-\xi t\right)
^{2}\right)  \label{16}
\end{equation}%
where $\rho =264$ $m$ is the standard deviation and the initial distribution
is $x_{0}=2$ $km$ away from the start. The initial concentration can be
obtained from (\ref{16}) by taking $t=0.$ At the boundaries the following
conditions are taken:%
\begin{equation*}
u(0,t)=u(L,t)=0
\end{equation*}%
where $L=9$ $km$. Since the velocity is $0.5$ $m/s$, the initial
distribution is transported $4.8$ $km$ after $9600$ $s$. Fig. \ref{fig1} shows this
transportation.%

\begin{figure}[ht]
\centering\includegraphics[scale=0.6]{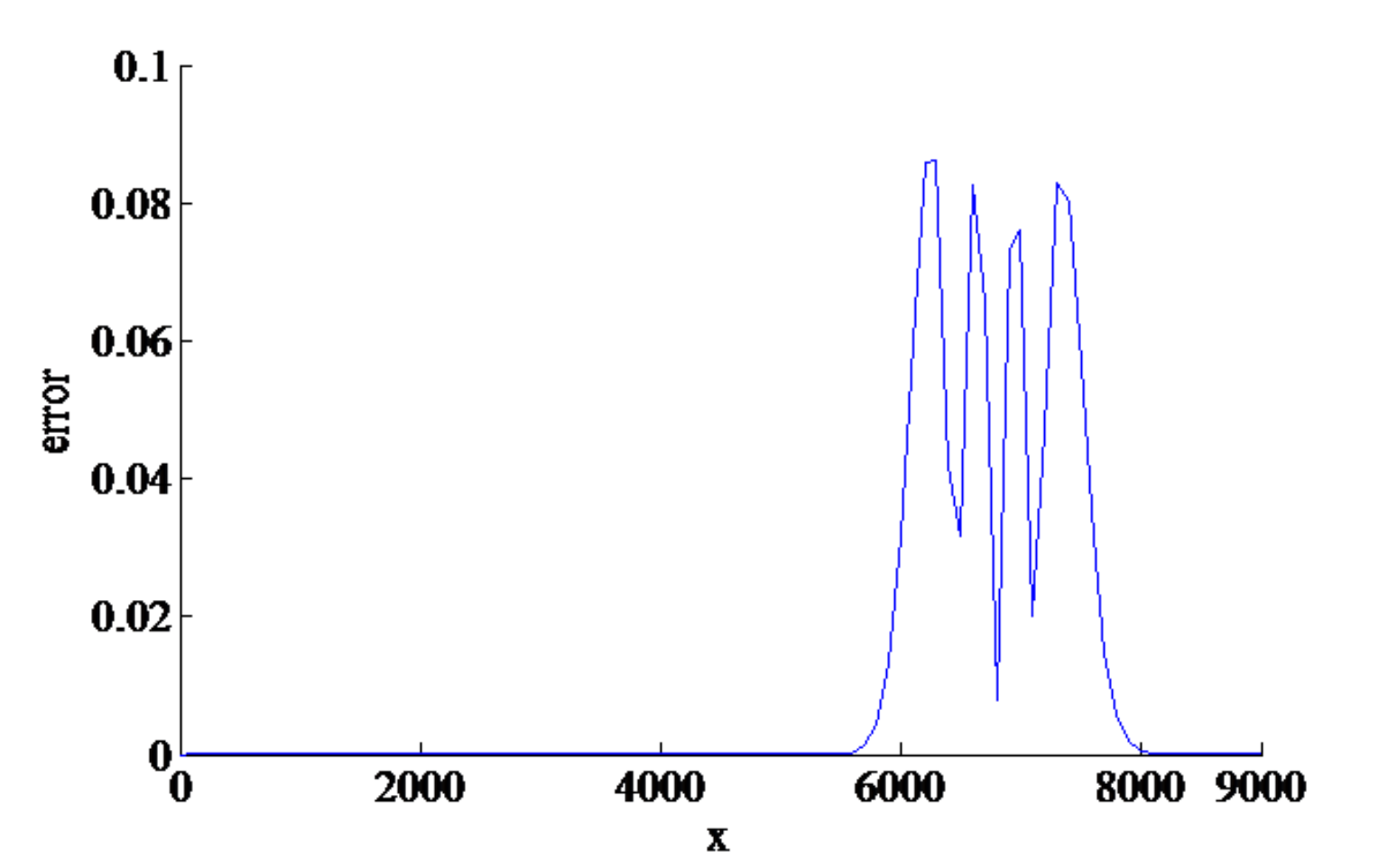}
\caption{{}{\protect \ Absolute error distributions at $t=9600$ with $C_{r}=0.25$ and $\Delta t=50$}}
\label{fig2}
\end{figure}

The obtained peak concentrations at $t=9600$ $s$ for various Courant numbers
are given in Table 2. As it is seen from the table that during the running
of time period, the calculations show that the results of EBSGM\textbf{\ }%
are considerably near to the exact value and in general, it leads to more
accurate results than the other methods. To see the errors along the whole
domain for various Courant numbers, Table 3 is documented. According to this
table, the results of EBSGM are mostly more accurate than the those produce
by the least square finite element method and extended cubic B-spline
collocation method for various Courant numbers. The absolute error
distributions of the EBSGM\textbf{\ }at $t=9600$ is illustrated in Fig. \ref{fig2}.
Maximum error occurred around the peak concentration.
\bigskip
\begin{equation*}
\begin{tabular}{cccccccc}
\multicolumn{8}{c}{Table 2: Peak concentrations at $t=9600$ $s$ for various $%
C_{r}$ and $\Delta t=50$} \\ \hline \hline
$C_{r}$ & $h$ & $p$ & EBSGM & \cite{szym} & \cite{dag} & \cite{gardner} &
Exact \\ \hline
0.25 & 100 & \multicolumn{1}{r}{6.8E-6} & 9.992 & 9.816 & 9.926 & 9.986 &
10.000 \\
0.50 & 50 & \multicolumn{1}{r}{13.6E-6} & 9.992 & 9.836 & 9.932 & 9.986 &
10.000 \\
0.75 & 33.3 & \multicolumn{1}{r}{2.04E-5} & 9.992 & 9.934 & 9.949 & 9.993 &
10.000 \\
1.00 & 25 & 3.59E-5 & 9.992 & 10.000 & 9.961 & 9.986 & 10.000 \\
1.50 & 16.6 & 4.91E-5 & 9.992 & 9.941 & 9.959 & 9.994 & 10.000 \\
2.00 & 12.5 & 7.18E-5 & 9.992 & 10.000 & 9.961 & 9.986 & 10.000 \\
3.20 & 7.8 & 7.50E-6 & 9.993 & 9.988 & 9.962 & 9.999 & 10.000 \\ \hline \hline
\end{tabular}%
\end{equation*}%
\bigskip
\begin{equation*}
\begin{tabular}{ccccccc}
\multicolumn{7}{c}{Table 3: Errors at $t=9600$ $s$ with $\xi =0.5$ $m/s$.}
\\ \hline \hline
$C_{r}$ & $h$ & $\Delta t$ & $p$ & EBSGM & \cite{dirk} & \cite{dag} \\ \hline
0.125 & 200 & 50 & \multicolumn{1}{r}{3.30E-6} & 1.63E-1 &
\multicolumn{1}{l}{1.29} & 5.18E-1 \\
0.25 & 100 & 50 & \multicolumn{1}{r}{6.80E-6} & 8.60E-2 & 3.25E-1 & 3.76E-1
\\
0.50 & 50 & 50 & \multicolumn{1}{r}{13.6E-6} & 9.07E-2 & 1.98E-1 & 3.73E-1
\\
0.50 & 10 & 10 & 1.53E-4 & 3.51E-3 & 7.51E-3 &  \\
0.50 & 1 & 1 & 3.04E-4 & 3.53E-5 & 7.50E-5 &  \\
0.50 & 0.5 & 0.5 & 3.40E-3 & 1.20E-5 & 1.88E-5 &  \\
0.75 & 33.3 & 50 & 2.04E-5 & 9.03E-2 &  & 3.76E-1 \\
1.00 & 25 & 50 & 3.59E-5 & 9.02E-2 &  & 3.79E-1 \\
1.50 & 16.6 & 50 & 4.91E-5 & 8.96E-2 &  & 3.78E-1 \\
2.00 & 12.5 & 50 & 7.18E-5 & 9.02E-2 &  & 3.79E-1 \\
3.20 & 7.8 & 50 & 7.50E-6 & 8.90E-2 &  & 3.80E-1 \\ \hline \hline
\end{tabular}%
\end{equation*}%
\bigskip \

\subsection{The distribution of an initial Gaussian pulse}

As a second test problem, we deal with both advection and diffusion. The
analytical solution to the one-dimensional ADE of a Gaussian pulse of unit
height over the domain $\left[ 0,9\right] $ is given as
\begin{equation}
u(x,t)=\frac{1}{\sqrt{4t+1}}\exp \left( -\frac{\left( x-x_{0}-\xi t\right)
^{2}}{\lambda \sqrt{4t+1}}\right)  \label{17}
\end{equation}%
where $\xi $ is the velocity, $\lambda $ is diffusion coefficient and $x_{0}$
is the centre of the initial Gaussian pulse \cite{sanka}.

\begin{figure}[ht]
\centering\includegraphics[scale=0.6]{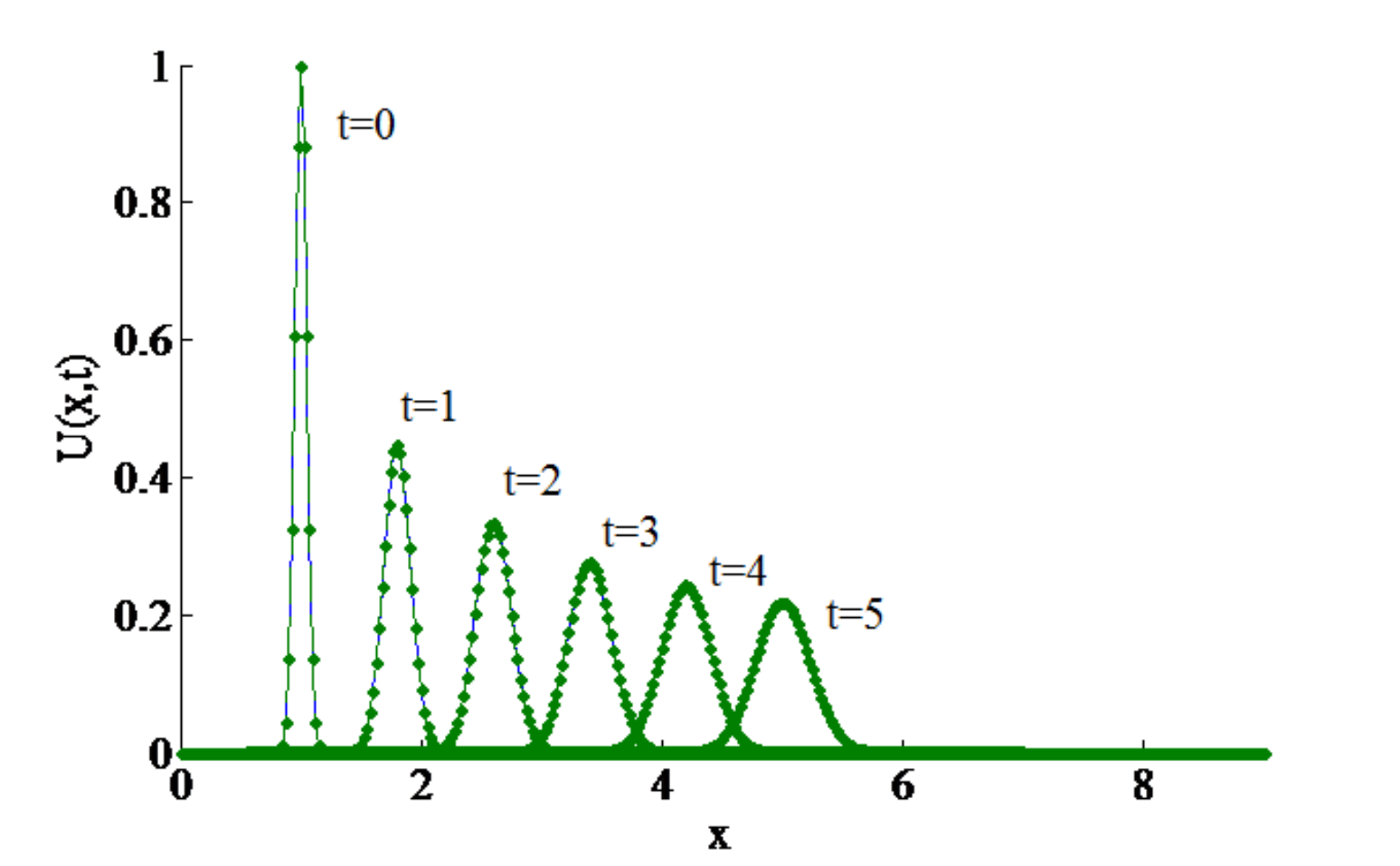}
\caption{{}{\protect \ Distribution of an initial Gaussian pulse}}
\label{fig3}
\end{figure}

\begin{figure}[ht]
\centering\includegraphics[scale=0.6]{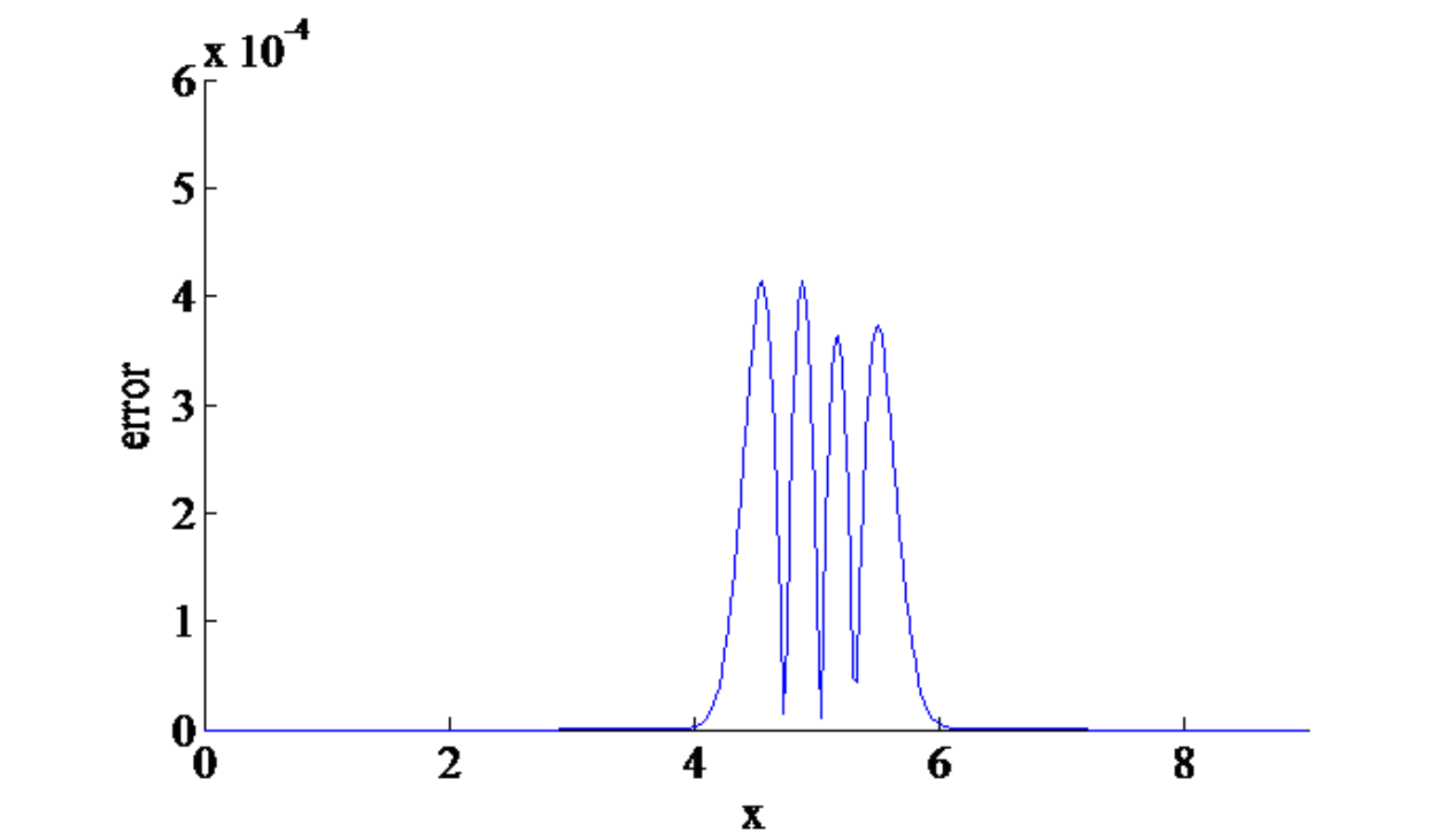}
\caption{{}{\protect \ Absolute error distributions at $t=5$ with $h=0.025, \Delta t=0.0125$}}
\label{fig4}
\end{figure}

The initial condition is chosen as the analytical value of the Eq. (\ref{17}%
) for $t=0$ and the boundary conditions are chosen as%
\begin{equation*}
u(0,t)=u(9,t)=0.
\end{equation*}%
The results presented here are computed for time step $\Delta t=0.0125$ $s$.
Parameters in the equation are used as $\lambda =0.005$ $m^{2}/s$ and $\xi
=0.8$ $m/s$. Fig.\ref{fig3} shows the behavior of the numerical and analytical
solutions (which are graphed with continuous lines) for various times until
the simulation terminating time $t=5$. Thus, the decay in time of the
initial pulse is modeled. So that the effect of the diffusion term has been
observed in this test problem. The absolute error distributions of the EBSGM%
\textbf{\ }at $t=5$ is illustrated in Fig. \ref{fig4}.%

For comparison, the advection-diffusion equation is solved for various
Courant numbers and computed errors at $t=5$ $s$ are presented in Table 4.%
\bigskip
\begin{equation*}
\begin{tabular}{ccccc}
\multicolumn{5}{c}{Table 4: Error norm at $t=5$, $\xi =0.8$ $m/s$, $\lambda
=0.005$ $m^{2}/s$, $\Delta t=0.0125$} \\ \hline \hline
$C_{r}$ & $h$ & EBSGM ($p=0.05286$) & Method I \cite{korkmaz} & Method II
\cite{korkmaz} \\ \hline
0.05 & \multicolumn{1}{l}{0.2} & 0.1326154 & 0.1253926 & 0.1361437 \\
0.10 & \multicolumn{1}{l}{0.1} & 0.0042464 & 0.0069553 & 0.0145554 \\
0.20 & \multicolumn{1}{l}{0.05} & 0.0008333 & 0.0012117 & 0.0002886 \\
0.40 & \multicolumn{1}{l}{0.025} & 0.0004134 & 0.0003071 & 0.0000181 \\
\hline \hline
\end{tabular}%
\end{equation*}

\section{Conclusion}

In this paper, we have proposed a new algorithm for the numerical solution
of the ADE. This algorithm is obtained by employing exponential B-spline
functions to the well known Galerkin finite element method. To see
achievement of the method is studied two test problems. The resulting
numerical solutions for various Courant numbers are compared with the
previous studies in Tables 3 and 4. Accordingly, we can say that the
proposed method give acceptable results.

\bigskip

\noindent \textbf{Acknowledgements}

\noindent The authors are grateful to The Scientific and Technological
Research Council of Turkey for financial support for their project 113F394
and granting scholarship for the author Melis' Ph.D. studies.

\end{document}